\documentclass[12pt]{article}
\usepackage{amsfonts,latexsym,amssymb}
 \newcommand{\beqn}{\begin{eqnarray}}
 \newcommand{\eeqn}{\end{eqnarray}}
 \newcommand{\be}{\begin{equation}}
 \newcommand{\ee}{\end{equation}}
 \newcommand{\ba}{\begin{array}}
 \newcommand{\ea}{\end{array}}
 \newcommand{\pa}{\partial}
 
 \newcommand{\ci}{\cite}
 \newcommand{\la}{\label}

\newcommand{\na}{\nabla}

 \newcommand{\De}{\Delta}

\def\R{{\rm I\kern-.1567em R}}
\def\M{{\rm I\kern-.1567em M}}

\def\div {{\rm div}}

\def\ess{{\rm ess}}

 \newtheorem{theorem}{Theorem}[section]
 
 \newtheorem{definition}[theorem]{Definition}
 
 \newtheorem{lemma}[theorem]{Lemma}
 
 \newtheorem{remark}[theorem]{Remark}

\begin{document}

\begin{center} {\bf Estimates of suitable weak solutions
   to the Navier-Stokes equations in critical Morrey spaces}\\
  \vspace{1cm}
 {\large
  G. Seregin}

 \end{center}

 \vspace{1cm}
 \noindent
 {\bf Abstract }
 We prove some estimates for suitable weak solutions to
 the nonstationary three-dimensional Navier-Stokes equations under assumptions
 that certain invariant functionals of the velocity field are bounded.

 \vspace {1cm}

\noindent {\bf 1991 Mathematical subject classification (Amer.
Math. Soc.)}: 35K, 76D.


\setcounter{equation}{0}
\section{Introduction  }

Consider the nonstationary 3D Navier-Stokes equations
 \be\la{i1}
\begin{array}{c}\pa_t v+v\cdot\na v-\De v
=-\nabla p,  \qquad
            \div \,v = 0\end{array}\ee
in the unit space-time cylinder $Q=B\times ]-1,0[\subset  \mathbb
R^3\times \mathbb R^1$. Here, $B(r)$ is the  ball of radius $r$ in
$\mathbb R^3$ centered at the space origin $x=0$, $Q(r)=B(r)\times
]-r^2,0[$ is a standard parabolic cylinder, $B=B(1)$, $Q=Q(1)$, $v$
and $p$ stand for the velocity and for the pressure, respectively.

It is  known that equations (\ref{i1}) are invariant with respect to
the following scaling (we call it the natural scaling)
$$v^\lambda(x,t)=\lambda v(\lambda x,\lambda^2t),\qquad
p^\lambda(x,t)=\lambda^2 p(\lambda x,\lambda^2t).$$

In the so-called $\varepsilon$-regularity theory, the important
role plays certain critical Morrey spaces. Their norms are
generated by functionals which are invariant with respect to the
natural scaling. Among such functionals, there are
$$C(r)=\frac 1{r^2}\int\limits_{Q(r)}|v|^3dz,\qquad A(r)=\ess
\sup\limits_{-r^2<t<0}\frac 1r\int\limits_{B(r)}|v(x,t)|^2dx,$$
$$E(r)=\frac 1r\int\limits_{Q(r)}|\na v|^2dz,\qquad H(r)=\frac 1{r^3}
\int\limits_{Q(r)}|v|^2dz,$$$$ D(r)=\frac
1{r^2}\int\limits_{Q(r)}|p|^\frac 32dz,\qquad  D_0(r)=\frac
1{r^2}\int\limits_{Q(r)}|p-[p]_{B(r)}|^\frac 32dz,$$
$$D_1(r)=\frac 1{r^\frac 32}\int\limits_{-r^2}^0
\Big(\int\limits_{B(r)}|\na p|^\frac 98dx\Big)^\frac 43dt,$$ where
$z=(x,t)$ is a point in space-time and
$$[p]_{B(r)}(t)=\frac 1{|B(r)|}\int\limits_{B(r)}p(x,t)dx.$$

All conditions of $\varepsilon$-regularity for the so-called
suitable weak solutions are stated with the help of those
functionals. For example, the famous Caffarelli-Kohn-Nirenberg
condition, see \ci{CKN},  reads as follows.
\begin{theorem}\la {it1}There is a universal positive constant
$\varepsilon$ with the following property. Assume that the pair
$v$ and $p$ is a suitable weak solution to the Navier-Stokes
equations in $Q$. If
\begin{equation}\label{i2}
    \sup\limits_{0<r\leq 1}E(r)\leq \varepsilon,
\end{equation}
then the space-time origin $z=0$ is a regular point of $v$.
\end{theorem}
Let us recall to the reader  definitions of suitable weak solutions
and regular points.

%As usual in local considerations, one should explain what kind of
%solutions to (\ref{i1}) is going to be treated. The most
%reasonable object is
%The  suitable weak solutions were introduced in \ci{CKN}, see also
%\ci{Sc1}, \ci{Li}, and \ci{LS}. We use more convenient definition
%given in \ci{Li}, see also \ci{LS}. Here, it is:
\begin{definition}\la{id2}
The pair $v$ and $p$ is called a suitable weak solution to the
Navier-Stokes equations in $Q$ if \be\la{i3}
    v\in L_{2,\infty}(Q)\cap W^{1,0}_2(Q),\qquad p\in L_\frac
    32(Q);
\ee \be\la{i4} \mbox{the Navier-Stokes equations hold in $Q$ in the
sense of distributions};\ee

\noindent for a.a. $t\in ]-1,0[$, the local energy inequality
$$\int\limits_B\varphi(x,t)|v(x,t)|^2dx+2\int\limits^t_{-1}
\int\limits_B\varphi|\na v|^2dxdt'$$
\begin{equation}\label{i5}
    \leq\int\limits^t_{-1}
\int\limits_B\Big\{|v|^2(\De
\varphi+\pa_t\varphi)+v\cdot\na\varphi(|v|^2+2p)\Big\}dxdt'
\end{equation}
holds for any non-negative test function $\varphi\in
C^\infty_0(\mathbb R^3\times \mathbb R^1)$ vanishing in a
neighborhood of the parabolic boundary of $Q$.
\end{definition}
\begin{definition}\la{id3} The point $z=0$ is called a regular point
of $v$ if there is a number $r\in ]0,1]$ such that $v$ is a H\"older
continuous function in $\overline{Q}(r)$.\end{definition} Here, the
following abbreviations are used:
$$L_{2,\infty}(Q)=L_\infty(-1,0;L_2(B)), \qquad
W^{1,0}_2(Q)=L_2(-1,0;W^1_2(B)),$$  $L_2(B)$ and $W^1_2(B)$ are
the usual Lebesgue and Sobolev spaces, respectively.
\begin{remark}\la{ir4}  Our definition of suitable
weak solutions belongs to F.-H. Lin \ci{Li}. It differs from more
general definition, given  by Caffarelly-Kohn-Niren- berg in
\ci{CKN}, by the very concrete choice of the space for the pressure.
To our opinion, such a choice seems to be more convenient to
treat.\end{remark}
\begin{remark}\la{ir5} Definition \ref{id3} of regular points is due
Ladyzhenskaya-Seregin \ci{LS}. In the most popular definition by
 Caffarelli-Kohn-Nirenberg, the H\"older space is replaced
 with the space of essentially bounded functions.\end{remark}

Roughly speaking, Theorem \ref{it1} and other similar statements say
that smallness of functionals, which are invariant with respect to
the natural scaling, is a sufficient condition for regularity.
Obviously, the next problem is to figure out what happens if  above
functionals are  bounded but not small.  This seems to be a subtle
and completely open question. However, there is one case, where the
answer is known and positive. It is the marginal case of the
so-called Ladyzhenskaya-Prodi-Serrin condition. Indeed, in the
Ladyzhenskaya-Prodi-Serrin condition, the key role plays the
functional $\|\cdot\|_{s,l,Q} $, which is the norm  of the mixed
Lebesgue space $L_{s,l}(Q)=L_l(-1,0;L_s(B))$. This norm is invariant
with the respect to the natural scaling if $3/s+2/l=1$ and $s\geq
3$. The regular case $s>3$ can be reduced to
 the smallness of the norm $\|v\|_{s,l,Q}$ with the help of
 the natural scaling and absolute continuity of Lebesgue's integral. So,
 the only case,
 which seems to be not reducible to the $\varepsilon$-regularity
 theory is $s=3$ and $l=+\infty$. It should be noticed that to
 treat $L_{3,\infty}$-case we had to develop a new method based on the
 unique continuation theory for parabolic equations, see \ci{ESS4}.

The aim of this paper is to contribute somehow to analysis of
smoothness of suitable weak solutions under additional assumptions
that certain functionals invariant with respect to natural scaling
are bounded. We hope that our results can be regarded as a
starting point for that analysis. Let us formulate them.
\begin{lemma}\la{il6} Assume that we are given a suitable weak solution
$v$ and $p$ in $Q$. Let, in addition,
\begin{equation}\label{i6}
    \sup_{0<r\leq 1}E(r)=E_0<+\infty.
\end{equation}
Then, there is a positive constant $d$ depending on $E_0$ only
such that\begin{equation}\label{i7}
    A^\frac 32(r)+C(r)+D_0^2(r)\leq d\Big(r^\frac 12(A^\frac 32(1)
    +D_0^2(1))+1\Big)
\end{equation}
for all $0<r\leq 1/4$.\end{lemma}

\begin{lemma} \la{il7} Suppose that the pair $v$ and $p$ is a suitable weak
 solution in $Q$. Let
\begin{equation}\label{i8}
 \sup_{0<r\leq 1}C(r)=C_0<+\infty.
\end{equation}
Then
\begin{equation}\label{i9}
    A(r)+D_0(r)+E(r)\leq c\Big(r^2 D_0(1)+C_0+C_0^\frac 23\Big)
\end{equation}
for all $0<r\leq 1/2$.
\end{lemma}
 Here and in what follows, $c$ is a positive universal constant.

\begin{lemma} \la{il8} Suppose that the pair $v$ and $p$ is a suitable weak
 solution in $Q$. Let
\begin{equation}\label{i10}
 \sup_{0<r\leq 1}A(r)=A_0<+\infty.
\end{equation}
Then there is a positive constant $e$ depending on $A_0$ only such
that
\begin{equation}\label{i11}
    C^\frac 43(r)+D_0(r)+E(r)\leq e\Big(r^2(D_0(1)+E(1))+1\Big)
\end{equation}
for all $0<r\leq 1/2$.
\end{lemma}
Statements similar to Lemmata \ref{il6}--\ref{il8} are proved by
Choe-Lewis in \ci{ChLe}, see Lemma 1 there. Our proof is different
and estimates are sharper.

  \noindent\textbf{Acknowledgement} The work was supported by the
Alexander von Humboldt Foundation, by  the RFFI grant 05-01-00941-a,
and by the CRDF grant RU-M1-2596-ST-04.

\setcounter{equation}{0}
\section{Preliminary inequalities }
There are three basic inequalities and their modifications. The
first of them is but a multiplicative inequality and has the
form\begin{equation}\label{p1}
    C(r)\leq c\Big[\Big(\frac \varrho r\Big)^3A^\frac 34(\varrho)
    E^\frac 34(\varrho)+\Big(\frac r\varrho\Big)^3A^\frac
    32(\varrho)\Big]
\end{equation}
for all $0<r\leq \varrho\leq 1$. The reader can find a proof of it
in \ci{LS}, see also \ci{Li}.

The second group of inequalities is a consequence of local energy
inequality (\ref{i5})
\begin{equation}\label{p2}
    A(R/2)+E(R/2)\leq c\Big[C^\frac 23(R)+C(R)+C^\frac
    13(R)D_0^\frac 23(R)\Big]
\end{equation}
for all $0<R\leq 1$. It follows from (\ref{i5}) directly. Another
version of the local energy inequality is demonstrated in \ci{LS}
\be\la {p3}A(R/2)+E(R/2)\leq c\Big[C^\frac 23(R)+C^\frac
    13(R)D_0^\frac 23(R)+A^\frac 12(R)C^\frac 13(R)E^\frac
    12(R)\Big]\ee
for all $0<R\leq 1$.

A kind of a decay estimate for the pressure is the third inequality.
There are a several versions of such decay estimate. One of the is
proved  in \ci{S6} and reads
\begin{equation}\label{p4}
    D_0(r)\leq c\Big[\Big(\frac r\varrho\Big)^\frac 52D_0(\varrho)+
    \Big(\frac \varrho r\Big)^2C(\varrho)\Big]
\end{equation}
for any $0<r\leq \varrho\leq 1$. However, in a number of cases, it
is more convenient to use a slightly different versions
\begin{equation}\label{p5}
 D_0(r)\leq c\Big[\Big(\frac r\varrho\Big)^\frac 52D_0(\varrho)+
    \Big(\frac \varrho r\Big)^2A^\frac 12(\varrho)E(\varrho)\Big]
\end{equation}
or
\begin{equation}\label{p6}
D_0(r)\leq c\Big[\Big(\frac r\varrho\Big)^\frac 52D_0(\varrho)+
    \Big(\frac \varrho r\Big)^3A^\frac 34(\varrho)
    E^\frac 34(\varrho)\Big].
\end{equation}
Both are valid for the same $r$ and $\varrho$ as in (\ref{p4}).

Inequalities (\ref{p5}) and (\ref{p6}) can be proved more or less in
the same way. To show the basic arguments, let us prove the first of
them. To this end, we decompose the pressure
\begin{equation}\label{p7}
    p=p_1+p_2
\end{equation}
in $B(\varrho)$ so that $p_1$ is a unique solution to the
variational identity
\begin{equation}\label{p8}
    \int\limits_{B(\varrho)}p_1\De \varphi dx =-
    \int\limits_{B(\varrho)}(\tau -\tau_\varrho):\na^2 \varphi dx,
\end{equation}
where $\varphi$ is an arbitrary test function from
$W^2_3(B(\varrho))$ satisfying the boundary condition $\varphi
|_{\pa B(\varrho)}=0$ and
$$\tau=(v-c_\varrho)\otimes (v-c_\varrho),\quad \tau_\varrho=
[(v-c_\varrho)\otimes (v-c_\varrho)]_{B(\varrho)}, \quad
c_\varrho=[v]_{B(\varrho)}.$$ Here, time $t$ is considered as a
parameter. Obviously, then,
\begin{equation}\label{p9}
    \De p_2=0
\end{equation}
in $B(\varrho)$.

We can easily find the estimate of $p_1$
$$ \int\limits_{B(\varrho)}|p_1|^\frac 32dx\leq c
 \int\limits_{B(\varrho)}|\tau -\tau_\varrho|^\frac 32dx.$$
 By the Galiardo-Nirenberg inequality,
$$\int\limits_{B(\varrho)}|p_1|^\frac 32dx\leq c
\Big(\int\limits_{B(\varrho)}|v-c_\varrho||\na v|dx\Big)^\frac 32$$
and thus
$$\int\limits_{B(\varrho)}|p_1|^\frac 32dx\leq c
\Big(\int\limits_{B(\varrho)}|v-c_\varrho|^2dx\Big)^\frac 34
\Big(\int\limits_{B(\varrho)}|\na v|^2dx\Big)^\frac 34.$$ On the
other hand, we can use the Poincar\'{e} inequality
$$\int\limits_{B(\varrho)}|v-c_\varrho|^2dx\leq c\varrho^2
\int\limits_{B(\varrho)}|\na v|^2dx$$ and the minimality property of
$c_\varrho$
$$\int\limits_{B(\varrho)}|v-c_\varrho|^2dx\leq
\int\limits_{B(\varrho)}|v|^2dx. $$ The latter relation leads to the
estimate
\begin{equation}\label{p10}
    \frac 1{\varrho^2}\int\limits_{-\varrho^2}^0
    \int\limits_{B(\varrho)}|p_1|^\frac 32dz\leq c E(\varrho)A^\frac
    12(\varrho).
\end{equation}

Since $p_2$ is a harmonic function in $B(\varrho)$, we have for
$0<r\leq \varrho /2$
$$\sup_{x\in B(r)}|p_2(x,t)-[p_2]_{B(r)}(t)
|^\frac 32\leq cr^\frac 32 \sup_{x\in B(\varrho /2)}|\na
p_2(x,t)|^\frac 32$$
\begin{equation}\label{p11}
   \leq c\Big(\frac r{\varrho^4}\int\limits_{B(\varrho)}
    |p_2(x,t)-[p_2]_{B(\varrho)}(t)|dx\Big)^\frac 32
\end{equation}
$$\leq \frac c{\varrho^3}\Big(\frac r\varrho\Big)^\frac 32
\int\limits_{B(\varrho)}
    |p_2(x,t)-[p_2]_{B(\varrho)}(t)|^\frac 32dx.$$

Next, by (\ref{p7}) and (\ref{p11}),
$$D(r)\leq \frac c{r^2}\int\limits_{Q(r)}|p_1-[p_1]_{B(r)}|^\frac 32dz
+\frac c{r^2}\int\limits_{Q(r)}|p_2-[p_2]_{B(r)}|^\frac 32dz$$
$$\leq \frac c{r^2}\int\limits_{Q(r)}|p_1|^\frac 32dz
+\frac c{r^2}\frac 1{\varrho^3}\Big(\frac r\varrho\Big)^\frac 32
\int\limits_{-r^2}^0r^3\int\limits_{B(\varrho)}
    |p_2(x,t)-[p_2]_{(\varrho)}(t)|^\frac 32dx$$
$$\leq c\Big(\frac \varrho r\Big)^2E(\varrho)A^\frac 12(\varrho)+
c\Big(\frac r\varrho\Big)^\frac 52\frac 1{\varrho^2}
\int\limits_{Q(\varrho)}|p_2-[p_2]_{B(\varrho)}|^\frac 32dz$$
$$\leq c\Big(\frac \varrho r\Big)^2E(\varrho)A^\frac 12(\varrho)+
c\Big(\frac r\varrho\Big)^\frac 52\Big[\frac 1{\varrho^2}
\int\limits_{Q(\varrho)}|p-[p]_{B(\varrho)}|^\frac 32dz$$$$ +\frac
1{\varrho^2} \int\limits_{Q(\varrho)}|p_1-[p_1]_{B(\varrho)}|^\frac
32dz\Big]$$
$$\leq c\Big[\Big(\frac r\varrho\Big)^\frac 52D_0(\varrho)
+\Big(\frac \varrho r\Big)^2E(\varrho)A^\frac 12(\varrho)\Big]$$ So,
inequality (\ref{p5}) is proved.

\setcounter{equation}{0}
\section{Proof of Lemma \ref{il6} }

So, assume that condition (\ref{i6}) holds. Then, as it follows from
(\ref{p1}), (\ref{p2}), and (\ref{i6}), we have
\begin{equation}\label{31}
    C(r)\leq c\Big[\Big(\frac\varrho r\Big)^3A^\frac
    34(\varrho)E_0^\frac 34+\Big(\frac r\varrho\Big)^3A^\frac
    32(\varrho)\Big]
\end{equation}
and
\begin{equation}\label{32}
    D_0(r)\leq c\Big[\Big(\frac r\varrho\Big)^\frac 52D_0(\varrho)
    +\Big(\frac\varrho r\Big)^2A^\frac 12(\varrho)E_0\Big].
\end{equation}

Introducing
$$\mathcal E(r)=A^\frac 32(r)+D_0^2(r),$$
we derive from local energy inequality (\ref{p3})
$$\mathcal E(r)\leq c\Big[C(2r)+C^\frac 12(2r)D_0(2r)+A^\frac
34(2r)C^\frac 12(2r)E_0^\frac 34\Big]+D^2_0(r)$$
\begin{equation}\label{33}
\leq c\Big[C(2r)+D_0^2(2r)+A^\frac 34(2r)C^\frac 12(2r)E_0^\frac
34\Big].
\end{equation}
Now, let us assume that $0<r\leq \varrho/2<\varrho\leq 1$.
Replacing $r$ with $2r$ in (\ref{31}) and (\ref{32}), we find from
(\ref{33})
$$\mathcal E(r)\leq c\Big[\Big(\frac\varrho r\Big)^3A^\frac
    34(\varrho)E_0^\frac 34+\Big(\frac r\varrho\Big)^3A^\frac
    32(\varrho)$$
$$+\Big(\frac r\varrho\Big)^5D_0^2(\varrho)
    +\Big(\frac\varrho r\Big)^4A(\varrho)E_0^2$$
$$+A^\frac 34(2r)\Big(\Big(\frac\varrho r\Big)^3A^\frac
    34(\varrho)E_0^\frac 34+\Big(\frac r\varrho\Big)^3A^\frac
    32(\varrho)\Big)^\frac 12E_0^\frac 34\Big]$$
$$\leq c\Big[\Big(\frac r\varrho\Big)^3A^\frac
    32(\varrho)+\Big(\frac r\varrho\Big)^5D_0^2(\varrho)
    +\Big(\frac r\varrho\Big)^\frac 32A^\frac 34(\varrho)E_0^\frac 34
    A^\frac 34(\varrho)\Big(\frac\varrho r\Big)^\frac 34$$
    $$+\Big(\frac\varrho r\Big)^{\frac 32+\frac 34}A^{\frac 34+\frac 38}(\varrho)
    E_0^{\frac 34+\frac 38}+\Big(\frac\varrho r\Big)^4A(\varrho)E_0^2
    +\Big(\frac\varrho r\Big)^3A^\frac
    34(\varrho)E_0^\frac 34\Big].$$
    Here, the obvious inequality $A(2r)\leq c\varrho A(\varrho)/r$
    has been used. Applying Young inequality with an arbitrary
    positive constant $\delta$, we show
    $$\mathcal E(r)\leq c\Big(\frac r\varrho\Big)^\frac 34(E_0^\frac 34+1)
    \mathcal E(\varrho)+c\delta \mathcal E(\varrho)$$
    $$+c(\delta)\Big(\Big(\frac\varrho r\Big)^6E_0^\frac 32+
    \Big(\frac\varrho r\Big)^{12}E_0^6+
    \Big(\frac\varrho r\Big)^{9}E_0^\frac 92\Big)\Big].$$
    Therefore,
    \begin{equation}\label{34}
\mathcal E(r)\leq c\Big[\Big(\frac r\varrho\Big)^\frac
34(E_0^\frac 34+1)+\delta\Big]
    \mathcal E(\varrho)+c(\delta)\Big(\frac\varrho
    r\Big)^{12}(E_0^6+E_0^\frac 92+E_0^\frac 32).
\end{equation}
Inequality (\ref{34}) holds for $r\leq \varrho/2$ and can be reduced
to  the form
  \begin{equation}\label{35}
\mathcal E(\vartheta\varrho)\leq c\Big[\vartheta^\frac
34(E_0^\frac 34+1)+\delta\Big]
    \mathcal E(\varrho)+c(\delta)\vartheta^{-12}(E_0^6+E_0^\frac 92+E_0^\frac 32)
\end{equation}
for any $0<\vartheta\leq 1/2$ and for any $0<\varrho\leq 1$.

Now, let us fix $\vartheta$ and $\delta$ in the following way
\begin{equation}\label{36}
    c\vartheta^\frac 14(E_0^\frac 34+1)<1/2,\quad 0<\vartheta\leq
    1/2,\quad c\delta<\vartheta^\frac 12/2.
\end{equation}
Obviously, $\vartheta$ and $\delta$ depend on $E_0$ only. So, we
have
\begin{equation}\label{37}
\mathcal E(\vartheta\varrho)\leq \vartheta^\frac 12 \mathcal
E(\varrho)+G
\end{equation}
for any $0<\varrho\leq 1$, where $\vartheta=\vartheta(E_0)$ and $G=
G(E_0)$.

Iterations of (\ref{37}) give us
$$\mathcal E(\vartheta^k\varrho)\leq\vartheta^\frac k2 \mathcal
E(\varrho)+cG$$ for any natural numbers $k$ and for any
$0<\varrho\leq 1$. Letting $\varrho=1$, we find
\begin{equation}\label{38}
\mathcal E(\vartheta^k)\leq\vartheta^\frac k2 \mathcal E(1)+cG
\end{equation}
for any natural numbers $k$.  It can be easily deduced from
(\ref{38}) that
\begin{equation}\label{39}
\mathcal E(r)\leq d_1(E_0)(r^\frac 12 \mathcal E(1)+1)
\end{equation}
for all $0<r\leq 1/2$. Now, for $C(r)$, we have from (\ref{31})
$$C(r)\leq c\Big[A^\frac 34(2r)E_0^\frac 34+A^\frac 32(2r)\Big]
\leq c\Big[A^\frac 32(2r)+E_0^\frac 32\Big]$$
$$\leq d_2(E_0)(\mathcal E(2r)+1)\leq d_3(E_0)(r^\frac 12\mathcal
E(1)+1).$$ So, Lemma \ref{il6} is proved.

\setcounter{equation}{0}
\section{Proof of Lemma \ref{il7} }

According to conditions (\ref{i8}) and inequality (\ref{p4}), the
following relation is valid:
$$D_0(r)\leq c\Big[\Big(\frac r\varrho\Big)^\frac
52D_0(\varrho)+\Big(\frac\varrho r\Big)^2C_0\Big]$$ for all
$0<r\leq\varrho\leq 1$. Letting $r=\vartheta\varrho$ with
$0<\vartheta\leq 1$, we find
$$D_0(\vartheta\varrho)\leq c\Big[\vartheta^\frac
52D_0(\varrho)+\vartheta^{-2}C_0\Big].$$ If we choose $\vartheta$ so
that $c\vartheta^\frac 12\leq 1$, then
$$D_0(\vartheta\varrho)\leq \vartheta^2D_0(\varrho)+cC_0.$$
Here, $c$ is a universal constant. After iterations, we arrive at
the inequality
$$D_0(\vartheta^k\varrho)\leq \vartheta^{2k}D_0(\varrho)
+cC_0$$ for any natural $k$. Setting $\varrho=1$, we find
$$D_0(\vartheta^k)\leq \vartheta^{2k}D_0(1)
+cC_0$$ for any natural $k$ or
\begin{equation}\label{41}
D_0(r)\leq cr^{2}D_0(1) +cC_0
\end{equation}
for any $0<r\leq 1$.

Next, by (\ref{p2}) and by (\ref{41}),
$$A(R/2)+E(R/2)\leq c\Big[C_0^\frac 23+C_0+C_0^\frac
13D_0^\frac 23(R)\Big]$$
\begin{equation}\label{42}
    \leq c\Big[D_0(R)+C_0^\frac 23+C_0\Big]
    \leq c\Big[R^2D_0(1)+C_0^\frac 23+C_0\Big]
\end{equation}
for all $0<R\leq 1$. So, we have
$$A(r)+D_0(r)+E(r)\leq c(r^2D_0(1)+C_0+C_0^\frac 23)$$
for all $0<r\leq 1/2$. Lemma \ref{il7} is proved.

\setcounter{equation}{0}
\section{Proof of Lemma \ref{il8} }

Here, we are going to use inequality (\ref{p1}) in the form
\begin{equation}\label{51}
    C(r)\leq c\Big[A^\frac 34(r)E^\frac 34(r)+A^\frac 32(r)\Big]
    \leq c\Big[A^\frac 34_0E^\frac 34(r)+A^\frac 32_0\Big].
\end{equation}
Next, according to local energy inequality (\ref{p2}), we have
$$\mathcal F(r)=E(r)+D_0(r)\leq c\Big[C^\frac 23(2r)+C(2r)
+D_0(2r)\Big]+D_0(r)$$
\begin{equation}\label{52}
    \leq c\Big[D_0(2r)+\Big(\frac \varrho r\Big)^2C(\varrho)
    +\Big(\frac \varrho r\Big)^\frac 43C^\frac 23(\varrho)\Big]
\end{equation}
for any $0<r\leq \varrho/2<\varrho\leq 1$. To prove (\ref{52}), the
inequality $C(r)\leq (\varrho/r)^2C(\varrho)$, $0<r\leq \varrho$,
has been used.

Now, by (\ref{p4}) and by (\ref{51}), (\ref{52}),
$$\mathcal F(r)\leq c\Big[\Big(\frac r\varrho\Big)^\frac
52D_0(\varrho)+\Big(\frac \varrho r\Big)^2C(\varrho)
    +\Big(\frac \varrho r\Big)^\frac 43C^\frac 23(\varrho)\Big]$$
$$\leq c\Big[\Big(\frac r\varrho\Big)^\frac
52D_0(\varrho)+\Big(\frac \varrho r\Big)^2\Big(A^\frac 34_0E^\frac
34(\varrho)+A^\frac 32_0\Big)$$$$+\Big(\frac \varrho r\Big)^\frac
43\Big(A^\frac 12_0E^\frac 12(\varrho)+A_0\Big)\Big].$$ Next, we
would like to exploit the fact that the power of $E(\varrho)$ is
less than one. To this end, the Young inequality with an arbitrary
positive constant $\delta$ is applied and we find
$$\mathcal F(r)\leq c\Big[\Big(\frac r\varrho\Big)^\frac
52D_0(\varrho)+\delta E(\varrho)\Big]$$$$+c(\delta)\Big[\Big(\frac
\varrho r\Big)^2A^\frac 32_0+\Big(\frac \varrho r\Big)^\frac 43A_0
+\Big(\frac \varrho r\Big)^8A^3_0+\Big(\frac \varrho r\Big)^\frac
83A_0\Big]$$
$$\leq c\Big(\Big(\frac r\varrho\Big)^\frac
52+\delta \Big)\mathcal F(\varrho)+c(\delta)\Big(\frac \varrho
r\Big)^8(A^3_0+A^\frac 32_0+A_0)$$ for any $0<r\leq
\varrho/2<\varrho\leq 1$. Now, let $r=\vartheta\varrho$ with
$0<\vartheta\leq 1/2$ and $0<\varrho\leq 1$.  As a result, we have
$$\mathcal F(\vartheta\varrho)\leq c(\vartheta^\frac 52+\delta)
\mathcal F(\varrho)+c(\delta)\vartheta^{-8}(A^3_0+A^\frac
32_0+A_0).$$ Fix $\vartheta$ and $\delta$ so that
$$c\vartheta^\frac 12\leq 1/2,\quad 0<\vartheta\leq 1/2,\quad
c\delta\leq \vartheta^2/2.$$ Clearly, $\vartheta$ and $\delta$ are
universal constants. This implies
$$\mathcal F(\vartheta\varrho)\leq \vartheta^2
\mathcal F(\varrho)+c(A^3_0+A^\frac 32_0+A_0).$$ Iterating the
latter relations and then letting $\varrho=1$, we arrive at the
estimate
$$\mathcal F(\vartheta^k)\leq \vartheta^{2k}
\mathcal F(1)+c(A^3_0+A^\frac 32_0+A_0)$$ being valid for any
natural number $k$ or
$$\mathcal F(r)\leq c r^{2}
\mathcal F(1)+c(A^3_0+A^\frac 32_0+A_0)$$ for any $0<r\leq 1/2$.

Next, we can derive from (\ref{51}) that
$$C^\frac 43(r)\leq c (A_0E(r)+A_0^2)\leq c (A_0\mathcal F(r)+A_0^2)$$
$$\leq e_1(A_0)(r^2\mathcal F(1)+1).$$
So, Lemma \ref{il8} is proved.

G. Seregin\\
Steklov Institute of Mathematics at St.Petersburg, \\
St.Peterburg, Russia


\begin{thebibliography}{99}

\bibitem {CKN}
Caffarelli, L., Kohn, R.-V., Nirenberg, L., Partial regularity of
suitable weak solutions of the Navier-Stokes equations, Comm. Pure
Appl. Math., Vol. XXXV (1982), pp. 771--831.

\bibitem {ChLe}
Choe, H. L., Lewis, J. L., On the singular set in the Navier-Stokes
equations, J. Functional Anal., 175(2000), pp. 348--369.

\bibitem{ESS4}
Escauriaza,L., Seregin, G.,  ~\v Sver\'ak, V.,.
$L_{3,\infty}$-Solutions to the Navier-Stokes equations and backward
uniqueness, Uspekhi Matematicheskih Nauk, v. 58, 2003, 2(350), pp. 3--44.
English translation in Russian Mathematical Surveys, 58(2003)2, pp.
211-250.


\bibitem {LS}
Ladyzhenskaya, O. A., Seregin, G. A., On partial regularity of
suitable weak solutions to the three-dimensional Navier-Stokes
equations, J.  math. fluid mech.,  1(1999), pp. 356-387.
\bibitem {Li}
Lin, F.-H., A new proof of the Caffarelly-Kohn-Nirenberg theorem,
Comm. Pure Appl. Math., 51(1998), no.3, pp. 241--257.
%\bibitem {Sc1}
%Scheffer, V., Hausdorff measure and the Navier-Stokes equations,
%Commun. Math. Phys., 55(1977), pp. 97--112.

\bibitem {S6}
 Seregin, G.A., On smoothness of $L_{3,\infty}$-solutions
 to the Navier-Stokes equations up to boundary,
 Mathematische Annalen,  332(2005), pp. 219-238.
%\bibitem{TX}
%Tian, G., Xin, Z., Gradient Estimation on Navier-Stokes Equations,
%Comm. Anal. Geom., 7(1999), No.2, 221--257.




\end{thebibliography}
\end{document}